\documentclass [11pt,reqno]{amsart}
\usepackage {amsmath,amssymb,epsfig,enumerate,verbatim}
\usepackage[all]{xy}

\newcommand{\C}{\mathbf{C}}

\newcommand{\R}{\mathbf{R}}

\renewcommand{\P}{\mathbf{P}}

\newcommand{\fX}{\mathfrak{X}}
\newcommand{\fY}{\mathfrak{Y}}

\newcommand{\cD}{\mathcal{D}}

\newcommand{\cH}{\mathcal{H}}

\newcommand{\Ltwo}{\mathrm{L}^2}

\renewcommand{\a}{\alpha}
\renewcommand{\b}{\beta}
\newcommand{\hodge}{H^{1,1}_\mathbf{R}}

\newcommand{\Nef}{\mathrm{Nef}}
\newcommand{\Psef}{\mathrm{Psef}}

\newcommand{\om}{\omega}

\newcommand{\cf}{{\rm cf.\ }} 
\newcommand{\eg}{{\rm e.g.\ }} 
\newcommand{\ie}{{\rm i.e.\ }}

\newcommand{\la}{\lambda}

\renewcommand{\=}{:=}

\DeclareMathOperator{\id}{id}

\DeclareMathOperator{\ord}{ord}

\numberwithin{equation}{section}       % Number formulas within sections

\newtheorem{prop} {Proposition} [section]
\newtheorem{thm}[prop] {Theorem} 
\newtheorem{defi}[prop] {Definition}
\newtheorem{lem}[prop] {Lemma}
\newtheorem{cor}[prop]{Corollary}
\newtheorem{prop-def}[prop]{Proposition-Definition}

\newtheorem*{mainthm}{Main Theorem} 
\theoremstyle{remark}

\newtheorem{rmk}[prop]{Remark}

\newtheorem*{ackn}{Acknowledgment} 

\title{Degree growth of meromorphic surface maps}

\date{\today}

\author{S{\'e}bastien Boucksom, Charles Favre, Mattias Jonsson}
\address{CNRS-Universit{\'e} Paris 7\\
  Institut de Math{\'e}matiques\\
  F-75251 Paris Cedex 05\\
  France\\
  and\\
  Graduate School of Mathematical Sciences\\
The University of Tokyo\\
3-8-1 Komaba, Meguro-ku, Tokyo, 153-8914\\
Japan} \email{boucksom@math.jussieu.fr}
\address{CNRS-Universit{\'e} Paris 7\\
  Institut de Math{\'e}matiques\\
  F-75251 Paris Cedex 05\\
  France\\
  and\\
 Unidade Mista CNRS-IMPA\\
 Estrada Dona Castorina 110\\
Rio de Janeiro / Brasil 22460-320}
\email{favre@math.jussieu.fr}
\address{Dept of Mathematics\\
  University of Michigan\\
  Ann Arbor, MI 48109-1109\\
  USA\\
and\\
Dept of Mathematics\\
KTH\\
SE-100 44 Stockholm\\
Sweden}

\email{mattiasj@umich.edu, mattiasj@kth.se} 
\thanks{First author
  supported by the Japanese Society for the Promotion of Science.
  Third author supported by the NSF, the Swedish Research Council and 
  the Gustafsson Foundation.}

\begin{document}

\begin{abstract}
  We study the degree growth of iterates of meromorphic selfmaps 
  of compact K\"ahler surfaces.
  Using cohomology classes on the Riemann-Zariski space
  we show that the degrees grow similarly to those of
  mappings that are algebraically stable on some bimeromorphic model.
\end{abstract}

\maketitle

%\setcounter{tocdepth}{1}
%\tableofcontents

%
%
%%%%%%%%%%%%%%%%%%%%%%%%%%%%%%%%%%%%%%%%%%%%%%%%%%%%%%%%%%%%%%%%%%%
%
%
\section*{Introduction}
Let $X$ be a compact K{\"a}hler surface and
$F:X\dashrightarrow X$ a dominant meromorphic mapping.
Fix a K\"ahler class $\om$ on $X$, normalized by
$(\om^2)_X=1$, and define the \emph{degree} of $F$ with respect 
to ${\om}$ to be the positive real number
\begin{equation*}
  \deg_{\om}(F)\=(F^*{\om}\cdot {\om})_X=(\om\cdot F_*\om)_X,
\end{equation*}
where $(\cdot)_X$ denotes the intersection form on $\hodge(X)$.
When $X=\P^2$ and ${\om}$ is the class of a line, this coincides with
the usual algebraic degree of $F$. 
%(In fact, not much is lost by restricting to this case.)
One can show that 
$\deg_{\om}(F^{n+m})\le 2\deg_{\om}(F^n)\deg_{\om}(F^m)$
for all $m,n$. Hence the limit 
\begin{equation*}
  \la_1\=\lim_{n\to\infty}\deg_{\om}(F^n)^{\frac1n},
\end{equation*}
exists. We refer to it as the \emph{asymptotic degree} 
of $F$.
It follows from standard arguments (see Proposition~\ref{prop:dinfty})
that $\la_1$ does not
depend on the choice of ${\om}$, that $\la_1$ is invariant under
bimeromorphic conjugacy, and that $\la_1^2\ge\la_2$,
where $\la_2$ is the topological degree of $F$.
\begin{mainthm}
  Assume that $\la_1^2>\la_2$. 
  Then there exists a constant $b=b(\om)>0$ such that
  \begin{equation*}
    \deg_{\om}(F^n)=b\la_1^n + O(\la_2^{n/2})
    \quad\text{as $n\to\infty$}.
  \end{equation*}
\end{mainthm}
The dependence of $b$ on $\omega$ can be made explicit: see
Remark~\ref{rem:bomega}.  For the polynomial map $F(x,y)=(x^d,x^dy^d)$
on $\C^2$ (with $\om$ the standard Fubini-Study form), one has
$\la_2=\la_1^2=d^2$, $\deg_\om(F^n) = nd^n$, hence the assertion in
the Main Theorem may fail when $\la_1^2=\la_2$.

Degree growth is an important component in the understanding of
the complexity and dynamical behavior of a selfmap and has been
studied in a large number of papers in both the mathematics
and physics literature. It is connected to topological entropy 
(see \eg~\cite{friedlandannals,guedjannals,guedjentropy,DSentropybound})
and controlling it is necessary in order to 
construct interesting invariant measures
and currents (see \eg~\cite{BF,FS2,RS,Sib}). 
Even in simple families of mappings, degree growth exhibits a rich
behavior: see \eg the papers by Bedford and Kim~\cite{BK1,BK2}, 
which also contain references to the physics literature.

In~\cite{FS2}, Forn{\ae}ss and Sibony connected the degree growth of
rational selfmaps to the interplay between contracted hypersurfaces
and indeterminacy points. In particular they proved that $\deg(F^n)$
is multiplicative iff $F$ is what is now often called (algebraically)
\emph{stable}. This analysis was extended to slightly more general
maps in~\cite{Nguyen}.  Forn{\ae}ss and Bonifant showed that
only countably many sequences $(\deg(F^n))_1^\infty$ can occur, but in
general the precise picture is unclear.

For bimeromorphic maps of surfaces, the situation is quite well
understood since the work of Diller and the second
author~\cite{DF}. Using the factorization into blowups and blowdowns,
they proved that any such map can be made stable by a bimeromorphic
change of coordinates. This reduces the study of degree growth to the
spectral properties of the induced map on the Dolbeault cohomology
$H^{1,1}$.  In particular it implies $\la_1$ is an algebraic integer,
that $\deg(F^n)$ satisfies an integral recursion formula and gives a
stronger version of our Main Theorem when $\la_1^2>1 (=\la_2)$.

In the case we consider, namely (noninvertible) meromorphic surface
maps, there are counterexamples to stability when
$\la_1^2=\la_2>1$~\cite{favremonomial}. 
It is an interesting 
(and probably difficult) 
question whether counterexamples also exist with 
$\la_1^2>\la_2>1$. 

Instead of looking for a particular birational model in which the
action of $F^n$ on $H^{1,1}$ can be controlled, we take a different
tack and study the action of $F$ on cohomology classes on all
modifications $\pi:X_\pi\to X$ at the same time.  This idea already
appeared in the study of cubic surfaces in~\cite{manin}, and was
recently used by Cantat as a key tool in his investigation of the
group of birational transformation of surfaces, see~\cite{cantat}.  In
the context of noninvertible maps, Hubbard and Papadopol~\cite{HP}
used similar ideas, but their methods apply only to a quite restricted
class of maps.

Here we show that $F$ acts (functorially) by pullback $F^*$ and
pushforward $F_*$ on the vector space $W\=\varprojlim\hodge(X_\pi)$
and on its dense subspace $C\=\varinjlim\hodge(X_\pi)$.  Compactness
properties of $W$ imply the existence of eigenvectors, having
eigenvalue $\la_1$ and certain positivity properties.

Following~\cite{DF} we then study the spectral properties of $F^*$ and
$F_*$ under the assumption $\la_1^2>\la_2$.  The space $W$ is too big
for this purpose, and we introduce a subspace $\Ltwo$ which is the
completion of $C$ with respect to the (indefinite) inner product
induced by the cup product, which is of Minkowski type by the Hodge
index theorem.  The Main Theorem then follows from the spectral
properties of $F^*$ and its adjoint $F_*$ on $\Ltwo$.

Using a different method, polynomial mappings of $\C^2$ were studied
in detail by the last two authors in~\cite{eigenval}: in that case
$\la_1$ is a quadratic integer.  However, our Main Theorem for
polynomial maps does not immediately follow from the analysis
in~\cite{eigenval}: the methods of the two papers can be viewed
as complementary.

The space $W$ above can be thought of as the Dolbeault 
cohomology $H^{1,1}$ of the \emph{Riemann-Zariski space} of $X$.
While we do not need the structure of the latter space in this paper, 
the general philosophy of considering all bimeromorphic models at the same 
time is very useful for handling asymptotic problems in 
geometry, analysis and dynamics: 
see~\cite{Diskant,cantat,manin} and~[FJ1-3]. 
In the present setting, it allows us to bypass the 
intricacies of indeterminacy points: 
heuristically, a meromorphic map becomes holomorphic on the 
Riemann-Zariski space.

The paper is organized in three sections. In the first we recall
some definitions and introduce cohomology classes on the 
Riemann-Zariski space. In the second, we study the actions of
meromorphic mappings on these classes. Finally, Section~3 deals
with the spectral properties of these actions under iteration,
concluding with the proof of the Main Theorem.

\smallskip
\noindent \emph{Remark on the setting}.  We chose to state our main
result in the context of a complex manifold, because the study of
degree growth is particularly important for applications to
holomorphic dynamics. 
However, our methods are purely algebraic so that our main
result actually holds in the case when $X$ is a projective surface
over any algebraically closed field of any characteristic, 
and $\om = c_1(L)$ for some ample line bundle. 
In this context, one has to replace
$\hodge (X)$ by the real N\'eron-Severi vector space, and work with the
suitable notion of pseudoeffective and nef classes, as defined
in~\cite[\S1.4, \S2.2]{Laz}.

\begin{ackn}
  We thank Serge Cantat and Jeff Diller for many useful remarks and
  the referees for a careful reading of the paper.
\end{ackn}
%
%
%%%%%%%%%%%%%%%%%%%%%%%%%%%%%%%%%%%%%%%%%%%%%%%%%%%%%%%%%%%%%%%%%%%
%
%
%\cleardoublepage
\section{Classes on the Riemann-Zariski space}
Let $X$ be a complex compact K{\"a}hler surface 
(for background see~\cite{BHPV}) and write
$\hodge(X)\=H^{1,1}(X)\cap H^2(X,\R)$.
%
%%%%%%%%%%%%%%%%%%%%%%%%%%%%%%%%%%%%%%%%%%%%%%%%%%%%%%%%%%%%%%%%%%%
%
\subsection{The Riemann-Zariski space}
By a \emph{blowup} of $X$, we mean a bimeromorphic morphism
$\pi:X_\pi\to X$ where $X_\pi$ is a smooth surface.  Up to
isomorphism, $\pi$ is then a finite composition of point blowups.  If
$\pi$ and $\pi'$ are two blowups of $X$, we say that $\pi'$
\emph{dominates} $\pi$ and write $\pi'\ge \pi$ if there exists a
bimeromorphic morphism $\mu:X_{\pi'}\to X_\pi$ such that
$\pi'=\pi\circ\mu$.  The \emph{Riemann-Zariski space} of $X$ is the
projective limit
\begin{equation*}
  \fX:=\varprojlim_\pi X_\pi.
\end{equation*}
While suggestive, the space $\fX$ is strictly speaking not
needed for our analysis and we refer 
to~\cite[Ch.VI, \S17]{ZS},~\cite[\S7]{Vaquie} for 
details on its structure.
%
%%%%%%%%%%%%%%%%%%%%%%%%%%%%%%%%%%%%%%%%%%%%%%%%%%%%%%%%%%%%%%%%%%%
%
\subsection{Weil and Cartier classes}
When one blowup $\pi'=\pi\circ\mu$ dominates another one $\pi$, we
have two induced linear maps $\mu_*:\hodge(X_{\pi'})\to\hodge(X_\pi)$
and $\mu^*:\hodge(X_\pi)\to\hodge(X_{\pi'})$, which satisfy the
projection formula $\mu_*\mu^*=\id$. This allows us to define the
following spaces.
\begin{defi}
  The space of Weil classes on $\fX$ is the projective limit
  \begin{equation*}
    W(\fX):=\varprojlim_\pi \hodge(X_\pi).
  \end{equation*}
  with respect to the push-forward arrows. The space of Cartier classes
  on $\fX$ is the inductive limit
  \begin{equation*}
    C(\fX):=\varinjlim_\pi\hodge(X_\pi).
  \end{equation*} 
  with respect to the pull-back arrows. 
\end{defi}
The space $W(\fX)$ is endowed with its projective limit
topology, \ie the coarsest topology for which the projection
maps $W(\fX)\to\hodge(X_\pi)$ are continuous.
There is also an inductive limit topology on
$C(\fX)$, but we will not use it. 

Concretely, a Weil class $\a\in W(\fX)$ is given by its
\emph{incarnations} $\a_\pi\in\hodge(X_\pi)$, compatible by
push-forward, that is $\mu_*a_{\pi'}=\a_\pi$ whenever
$\pi'=\pi\circ\mu$.  The topology on $W(\fX)$
is characterized as follows: 
a sequence (or net\footnote{A net is a family indexed by a directed set,
see~\cite{folland}.}) $\a_j\in W(\fX)$ converges to $\a\in W(\fX)$ iff
$\a_{j,\pi}\to\a_\pi$ in $\hodge(X_\pi)$ for each blowup $\pi$.

The projection formula recalled above shows that there is an injection
$C(\fX)\subset W(\fX)$, so that a Cartier class is in particular a
Weil class. In fact, if $\a\in\hodge(X_\pi)$ is a class in some
blow-up $X_\pi$ of $X$, then $\a$ defines a Cartier class, also
denoted $\a$, whose incarnation $\a_{\pi'}$ in any blowup
$\pi'=\pi\circ\mu$ dominating $\pi$ is given by $\a_{\pi'}=\mu^*\a$.
We say that $\a$ is \emph{determined} in $X_\pi$. (It is then also
determined in $X_{\pi'}$ for any blowup dominating $\pi$). Each
Cartier class is obtained that way. The space $C(\fX)$ is dense in
$W(\fX)$: if $\alpha$ is a given Weil class, the net $\a_\pi$ of
Cartier classes determined by the incarnations of $\a$ on all models
$X_\pi$ tautologically converges to $\a$ in $W(\fX)$.
\begin{rmk}
  The spaces of Weil classes and Cartier classes are 
  denoted by $Z_{\cdot}(X)$ and $Z^{\cdot}(X)$
  by Manin~\cite{manin}. He views these classes as living on the 
  ``bubble space'' $\varinjlim X_\pi$ rather than the
  Riemann-Zariski space $\varprojlim X_\pi$.
\end{rmk}

\subsection{Exceptional divisors}
This section can be skipped on a first reading, the main technical issue
being Proposition~\ref{p:conv}, to be used for the proof of
Theorem~\ref{thm:eigenclass}.

The spaces $C(\fX)$ and $W(\fX)$ are clearly bimeromorphic invariants
of $X$.  Once the model $X$ is fixed, an alternative and somewhat more
explicit description of these spaces can be given in terms of
exceptional divisors.
\begin{defi} 
  The set $\cD$ of 
  exceptional primes over $X$ is defined as the set of all exceptional
  prime divisors of all blow-ups $X_\pi\to X$ modulo the following
  equivalence relation: two divisors $E$, $E'$ on $X_\pi$ and $X_{\pi'}$
  are equivalent if the induced meromorphic map
  $X_\pi\dashrightarrow X_{\pi'}$ sends $E$ onto $E'$.
\end{defi}
When $X$ is a projective surface, $\cD$ is the 
set of divisorial valuations on the function field $\C(X)$
whose center on $X$ is a point.

If $E\in\cD$ is an exceptional prime and $X_\pi$ is any
model of $X$, one can consider the \emph{center} of $E$ on $X_\pi$,
denoted by $c_\pi(E)$. It is a subvariety defined as follows: choose a
blow-up $\pi'\ge\pi$ such that $E$ appears as a curve on
$X_{\pi'}$. Then $c_\pi(E)$ is defined as the image of $E\subset
X_{\pi'}$ by the map $X_{\pi'}\to X_{\pi}$. It does not depend on the
choice of $\pi'$, and is either a point or an irreducible curve. In
this 2-dimensional setting, there is a unique minimal blow-up $\pi_E$
such that $c_\pi(E)$ is a curve iff $\pi\ge\pi_E$ (in particular
$c_{\pi_E}(E)$ is a curve).

Using these facts, one can construct an explicit basis for the vector
space $C(\fX)$ as follows
(compare~\cite[Proposition~35.6]{manin}). Let $\a_E\in C(\fX)$ be the
Cartier class determined by the class of $E$ on $X_{\pi_E}$.  Write
$\R^{(\cD)}$ for the direct sum $\oplus_\cD \R$, or
equivalently for the space of real-valued functions on $\cD$ with
finite support.

\begin{prop}\label{p:direct} 
  The set $\{\a_E\ |\ E\in\cD\}$ is
  a basis for the vector space of Cartier classes $\a\in C(\fX)$ that
  are exceptional over $X$, \ie whose incarnations on $X$ vanish.
  In other words, the map 
  $\hodge(X)\oplus\R^{(\cD)}\to C(\fX)$ sending 
  $\alpha\in\hodge(X)$ to the Cartier class it
  determines and $E\in\cD$ to $\a_E$ is an isomorphism.
\end{prop}
%If $\pi$ is a given blow-up with exceptional primes $E_1,\dots,E_N$,
%the map $C(\fX)\to C(\fX)$ $\a\mapsto\a_\pi$ corresponds to the
%projection $\R^{(\cD)}\to\R^N\subset\R^{(\cD)}$ onto the
%$\a_{E_i}$-components.

We now describe $W(\fX)$ in terms of exceptional primes.  If $\a\in
W(\fX)$ is a given Weil class, let $\alpha_X\in\hodge(X)$ be its
incarnation on $X$. For each $\pi$, the Cartier class
$\alpha_\pi-\alpha_X$ is determined on $X_\pi$ by a unique
$\R$-divisor $Z_\pi$ exceptional over $X$. If $E$ is a
$\pi$-exceptional prime, we set $\ord_E(\a):=\ord_E(Z_\pi)$ so that
$Z_\pi= \sum_E \ord_E(Z_\pi) E$. It is easily seen to depend only on
the class of $E$ in $\cD$.  Let $\R^\cD$ denote the (product) space of
all real-valued functions on $\cD$. We obtain a map
$W(\fX)\to\hodge(X)\times\R^\cD$, which is easily seen to be a
bijection, and even naturally a homeomorphism as the following
straightforward lemma shows.
\begin{lem}\label{p:alter} 
  A net $\a_j\in W(\fX)$ 
  converges to $\a\in W(\fX)$ iff $\a_{j,X}$ converges to $\a_X$ in
  $\hodge(X)$ and $\ord_E(\a_j)\to\ord_E(\a)$ for each exceptional
  prime $E\in\cD$.
\end{lem}

A result of Zariski 
(\cf \cite[Theorem~3.17]{kollar},~\cite[Proposition~1.12]{valtree})
states that the process of successively blowing-up the center of a
given exceptional prime $E\in\cD$ starting from any given
model must stop after finitely many steps with the center becoming a
curve. In other words, if 
$X=X_0\leftarrow X_1\leftarrow X_2\leftarrow\dots$ 
is an infinite sequence of blow-ups such that the
center of each blow-up $X_n\leftarrow X_{n+1}$ meets $c_{X_n}(E)$,
then $X_n$ must dominate $X_{\pi_E}$ for $n$ large enough.
Using this
result, we record the following fact to be used later on:
\begin{prop}\label{p:conv} 
  Let $X=X_0\leftarrow X_1\leftarrow X_2\leftarrow\dots$ be an
  infinite sequence of blow-ups, and for each $n$ suppose $\a_n\in
  C(\fX)$ is a Cartier class which is determined in $X_{n+1}$ and
  whose incarnation on $X_n$ is zero.  Then $\a_n\to 0$ in $W(\fX)$ as
  $n\to\infty$.
\end{prop}
\begin{proof} 
  In view of Proposition~\ref{p:alter}, we have to show that for every
  given exceptional prime $E\in\cD$, $\ord_E(\a_n)$ converges
  to $0$ as $n\to\infty$. In fact, we claim that $\ord_E(\a_n)=0$ for
  $n\ge n(E)$ large enough. Indeed, according to Zariski's result,
  there are two possibilities: either there exists $N$ such that
  $c_{X_N}(E)$ is a curve, or there exists $N$ such that the center of
  the blow-up $X_{n+1}\to X_n$ does not meet $c_{X_n}(E)$ for all $n\ge
  N$. In the first case, it is clear that $\ord_E(\a_n)=0$ for $n\ge
  N$, since $\a_n$ is exceptional over $X_N$. In the second case, the
  center of $E$ on $X_n$ does not meet the exceptional divisor of
  $X_n\to X_{n-1}$ for $n>N$, which supports the exceptional class
  $\a_n$, thus $\ord_E(\a_n)=0$ for $n>N$ as well.
\end{proof} 
%
%%%%%%%%%%%%%%%%%%%%%%%%%%%%%%%%%%%%%%%%%%%%%%%%%%%%%%%%%%%%%%%%%%%
%
\subsection{Intersections and $\Ltwo$-classes}\label{L2}
For each $\pi$, the intersection pairing
$\hodge(X_\pi)\times\hodge(X_\pi)\to\R$ will be denoted by
$(\a\cdot\b)_{X_\pi}$. It is non-degenerate, and satisfies the
projection formula:
$(\mu_*\a\cdot\b)_{X_\pi}=(\a\cdot\mu^*\b)_{X_{\pi'}}$ if
$\pi'=\pi\circ\mu$. It thus induces a pairing 
$W(\fX)\times C(\fX)\to\R$ which will simply 
be denoted by $(\a\cdot\b)$.
\begin{prop}\label{p:dual}
  The intersection 
  pairing induces a topological isomorphism between $W(\fX)$ and
  $C(\fX)^*$ endowed with its weak-$*$ topology.
\end{prop}
\begin{proof} 
  A linear form $L$ on $C(\fX)=\varinjlim_\pi\hodge(X_\pi)$ 
  is the same thing as a collection of linear forms $L_\pi$ on
  $\hodge(X_\pi)$, compatible by restriction. Now such a collection is
  by definition an element of the projective limit
  $\varprojlim_\pi\hodge(X_\pi)^*$, which is identified to $W(\fX)$ via
  the intersection pairing. This shows that the intersection pairing
  identifies $W(\fX)$ with the dual of $C(\fX)$ endowed with its weak-$*$
  topology.
\end{proof}
The intersection pairing defined above restricts to a non-degenerate
quadratic form on $C(\fX)$, denoted by $\alpha\mapsto(\a^2)$. 
However, it does \emph{not} extend to a continuous
quadratic form on $W(\fX)$. For instance, if $z_1,z_2,...$ is a
sequence of distinct points on $X$ and $\pi_n$ denotes the blow-up of
$X$ at $z_1,...,z_n$, with exceptional divisor $F_n=E_1+...+E_n$, we
have $(F_n^2)=-n$, but $\{F_n\}\in C(\fX)$ converges in $W(\fX)$.  We
thus introduce the maximal space to which the intersection form
extends:
\begin{defi}
  The space 
  of $L^2$ classes $\Ltwo(\fX)$ is defined as the completion of $C(\fX)$
  with respect to the intersection form.
\end{defi}
The usual setting to perform a completion is that of a definite
quadratic form on a vector space, which is not the case of the
intersection form on $C(\fX)$. However, the Hodge index theorem
implies that it is of Minkowski type, and it is easy to show that the
completion exists in that setting.

Let us be more precise: if $\omega\in C(\fX)$ is a given class with
$(\omega^2)>0$, the intersection form is negative definite on its
orthogonal complement 
$\omega^{\perp}:=\{\a\in C(\fX)\ |\ (\a\cdot\omega)=0\}$ 
as a consequence of the Hodge index theorem
applied to each $\hodge(X_\pi)$. We have an orthogonal decomposition
$C(\fX)=\R\omega\oplus\omega^{\perp}$, and we then let
$\Ltwo(\fX):=\R\omega\oplus\overline{\omega^{\perp}}$, where
$\overline{\omega^{\perp}}$ is the completion in the usual sense of
$\omega^{\perp}$ endowed with the negative definite quadratic form
$(\alpha^2)$. Note that $t\omega\oplus\a\mapsto t^2-(\a^2)$ is then a
norm on $\Ltwo(\fX)$ that makes it a Hilbert space, but this norm
depends on the choice of $\omega$. However, the topological vector
space $\Ltwo(\fX)$ does not depend on the choice of $\omega$.

\label{blob}In fact, the completion can be characterized by the following
universal property: if $(Y,q)$ is a complete topological vector space
with a continuous non-degenerate quadratic form of Minkowski type, 
any isometry $T:C(\fX)\to Y$ continuously extends to $\Ltwo(\fX)\to Y$.

The intersection form on $\Ltwo(\fX)$ is also of Minkowski
type, so that it satisfies the Hodge index theorem: if a non-zero
class $\a\in\Ltwo(\fX)$ satisfies $(\a^2)>0$, then the intersection form
is negative definite on $\a^{\perp}\subset\Ltwo(\fX)$.
\begin{rmk}\label{Rem:21}
  The direct sum decomposition
  $C(\fX)=\hodge(X)\oplus\R^{(\cD)}$ of
  Proposition~\ref{p:direct} is orthogonal with respect to the
  intersection form. 
  Furthermore, the intersection form is negative definite on
  $\R^{(\cD)}$ and $\{\a_E\ |\ E\in\cD\}$ forms  an orthonormal
  basis for $-(\alpha^2)$. Indeed, the center of $E\in\cD$ on the
  minimal model $X_{\pi_E}$ on which it appears is necessarily the last
  exceptional divisor to have been created in any factorization of
  $\pi_E$ into a sequence of point blow-ups, thus it is a $(-1)$-curve.
  
  Using this, one sees that $\Ltwo(\fX)$ is isomorphic
  to the direct sum $\hodge(X)\oplus\ell^2(\cD)\subset W(\fX)$ 
  where $\ell^2(\cD)$ denotes the set of real-valued square-summable
  functions $E\mapsto a_E$ on $\cD$.
\end{rmk}

The different spaces we have introduced so far are related as follows.
\begin{prop}\label{p:charac} 
  There is a 
  natural continuous injection $\Ltwo(\fX)\to W(\fX)$, 
  and the topology
  on $\Ltwo(\fX)$ induced by the topology of $W(\fX)$ coincides
  with its weak topology as a Hilbert space.
  
  If $\a\in W(\fX)$ is a given Weil class, then the intersection number
  $(\a_\pi^2)$ is a decreasing function of $\pi$, and $\a\in\Ltwo(\fX)$
  iff $(\a_\pi^2)$ is bounded from below, in which case
  $(\a^2)=\lim_\pi(\a_\pi^2)$.
\end{prop}
\begin{proof}
  The injection 
  $\Ltwo(\fX)\to W(\fX)$ is dual to the dense injection
  $C(\fX)\subset\Ltwo(\fX)$. By Proposition~\ref{p:dual}, a net
  $\a_k\in\Ltwo(\fX)$ converges to $\a\in\Ltwo(\fX)$ in the topology
  induced by $W(\fX)$ iff $(\a_k\cdot\b)\to(\a\cdot\b)$ for each $\b\in
  C(\fX)$. Since $C(\fX)$ is dense in $\Ltwo(\fX)$, this
  implies $\alpha_k\to\alpha$ weakly in $\Ltwo(\fX)$.
  
  For the last part, one can proceed using the abstract definition
  of $\Ltwo(\fX)$ as a completion, but it is more transparent
  to use the explicit representation of Remark~\ref{Rem:21}.
  For any $\pi$, we have 
  $\a_\pi=\a_X+\sum_{E\in\cD_\pi}(\a\cdot\a_E)\a_E$, 
  where $\cD_\pi\subset\cD$
  is the set of exceptional primes of $\pi$.
  Then $(\a_\pi^2)=(\a_X^2)-\sum_{E\in\cD_\pi}(\a\cdot\a_E)^2$,
  which is decreasing in $\pi$.
  It is then clear that $\a\in\Ltwo(\fX)$ iff
  $(\a_\pi^2)$ is uniformly bounded from below 
  and $(\a^2)=\lim(\a_\pi^2)$.
\end{proof}
% 
%%%%%%%%%%%%%%%%%%%%%%%%%%%%%%%%%%%%%%%%%%%%%%%%%%%%%%%%%%%%%%%%%%%
%
\subsection{Positivity}
Recall that a class in $\hodge(X)$ is \emph{psef} (pseudoeffective)
if it is the class of a closed positive $(1,1)$-current on $X$.
It is \emph{nef} (numerically effective) if it is the limit of 
K\"ahler classes. Any nef class is psef.
The cone in $\hodge(X)$ consisting of psef classes is
strict: if $\a$ and $-\a$ are both psef, then $\a=0$.

If $\pi'=\pi\circ\mu$ is a blowup dominating some other blowup $\pi$,
then $\a\in\hodge(X_\pi)$ is psef (nef) iff $\mu^*\a\in\hodge(X_{\pi'})$
is psef (nef). On the other hand, if $\a'\in\hodge(X_{\pi'})$ 
is psef (nef), then so is $\mu_*\a'\in\hodge(X_\pi)$.
(For the nef part of the last assertion it is important that we 
work in dimension two.)
\begin{defi} 
  A Weil class $\a\in W(\fX)$ is 
  psef (nef) if its incarnation $\a_\pi\in\hodge(X_\pi)$ is 
  psef (nef) for any blowup $\pi:X_\pi\to X$.
\end{defi}
We denote by $\Nef(\fX)\subset\Psef(\fX)\subset W(\fX)$ the convex
cones of nef and psef classes.  The remarks above imply that a Cartier
class $\a\in C(\fX)$ is psef (nef) iff $\a_\pi\in\hodge(X_\pi)$ is
psef (nef) for one (or any) $X_\pi$ in which $\a$ is determined. We write
$\a\ge\b$ as a shorthand for $\a-\b\in W(\fX)$ being psef.
\begin{prop}\label{p:compact}
  The nef cone $\Nef(\fX)$ and the psef cone $\Psef(\fX)$ are 
  strict, closed, convex cones in $W(\fX)$ with compact bases.
\end{prop}
\begin{proof}
  The nef (resp.\ psef) cone is the projective limit of the
  nef (resp.~psef) cones of each $\hodge(X_\pi)$.
  These are strict, closed, convex cones with compact bases, 
  so the result follows from the Tychonoff theorem.
\end{proof}

Nef classes satisfy the following monotonicity property:
\begin{prop}\label{p:incr} 
  If $\a\in W(\fX)$ 
  is a nef Weil class, then $\a\le\a_\pi$ for each $\pi$. 
  In particular,
  $\a_\pi\neq 0$ for each $\pi$ unless $\a=0$.
\end{prop} 
\begin{proof} 
  By induction on the number
  of blow-ups, it suffices to prove that $\a_{\pi'}\leq\mu^*\a_\pi$
  when $\pi'=\pi\circ\mu$ and 
  $\mu$ is the blowup of a point in $X_\pi$. But then
  $\mu^*\a_\pi=\a_{\pi'}+cE$, where $E$ is the class of
  the exceptional divisor
  and $c=(\a_{\pi'}\cdot E)\ge0$. To get the second point, note that
  $\a_\pi=0$ for some $\pi$ implies $\a\leq 0$. On the other hand, 
  $\a\ge 0$ as $\a$ is nef. Since $\Psef(\fX)$ is a
  strict cone, we infer $\a=0$.
\end{proof}

\begin{prop}\label{p:monotone} 
  The nef
  cone $\Nef(\fX)$ is contained in $\Ltwo(\fX)$. If $\a_i\ge\b_i$,
  $i=1,2$ are nef classes, then we have
  $(\a_1\cdot\a_2)\ge(\b_1\cdot\b_2)\ge 0$.
\end{prop}
\begin{proof} 
  If $\a\in W(\fX)$ 
  is nef, each incarnation $\a_\pi$ is nef, and thus $(\a_\pi^2)\geq 0$,
  so that $\a\in\Ltwo(\fX)$ by Proposition~\ref{p:charac}, with
  $(\a^2)=\inf_\pi(\a_\pi^2)\geq 0$. To get the second point, note that
  $(\a_1\cdot\a_2)\geq (\a_1\cdot\b_2)$ since $\a_2-\b_2$ is psef and
  $\a_1$ is nef, and similarly $(\a_1\cdot\b_2)\ge(\b_1\cdot\b_2)$.
\end{proof}
 
These two propositions together show that if $\omega\in C(\fX)$ is a
Cartier class determined by a K\"ahler class down on $X$, then
$(\a\cdot\omega)>0$ for any non-zero nef class $\a\in W(\fX)$.

\begin{prop}\label{p:nefbounds}
  We have $2\,(\a\cdot\b)\, \a\ge(\a^2)\, \beta$ for any nef Weil
  classes $\a,\beta\in W(\fX)$.  In particular, if $\omega\in C(\fX)$
  is determined by a K\"ahler class on $X$ normalized by
  $(\omega^2)=1$, we have, for any non-zero nef Weil class $\a$:
  \begin{equation}\label{eq:nefbounds}
    \frac{(\a^2)}{2(\a\cdot {\om})}{\,\om}
    \le\a\le2(\a\cdot{\om})\,\om.
  \end{equation}
\end{prop}
\begin{proof}
  The second assertion is a special case of the first one. To prove
  the first one, we may assume $(\a\cdot\b)>0$, or else $\a$ and $\b$
  are proportional by the Hodge index theorem and the result is clear.
  It is a known fact (see the remark after Theorem~4.1
  in~\cite{BoucksomZariski}) that if $\gamma\in C(\fX)$ is a Cartier
  class with $(\gamma^2)\ge0$, then either $\gamma$ or $-\gamma$ is
  psef.  In view of Proposition~\ref{p:charac}, the same result is
  true for any $\gamma\in\Ltwo(\fX)$.  Apply this to $\gamma=\a-t\b$,
  where $t=\frac12(\a\cdot\a)/(\a\cdot\b)$.  As
  $(\gamma\cdot\gamma)\ge0$ and $(\gamma\cdot\a)\ge0$, 
  $\gamma$ must be psef.
\end{proof}
% 
%%%%%%%%%%%%%%%%%%%%%%%%%%%%%%%%%%%%%%%%%%%%%%%%%%%%%%%%%%%%%%%%%%%
%
\subsection{The canonical class}
The canonical class $K_\fX$ is the Weil class whose incarnation in any
blowup $X_\pi$ is the canonical class $K_{X_\pi}$. It is not Cartier
and does not even belong to $\Ltwo(\fX)$.  However, $K_{X_{\pi'}}\ge
K_{X_\pi}$ whenever $\pi'\ge \pi$, and $K_\fX$ is the smallest Weil
class dominating all the $K_{X_\pi}$. This allows us to intersect
$K_\fX$ with any nef Weil class $\a$ in a slightly ad-hoc way: we set
$(\a\cdot K_\fX)\=\sup_\pi(\a_\pi\cdot
K_{X_\pi})_{X_{\pi}}\in\R\cup\{+\infty\}$.
%
%
%%%%%%%%%%%%%%%%%%%%%%%%%%%%%%%%%%%%%%%%%%%%%%%%%%%%%%%%%%%%%%%%%%%
%
%
\section{Functorial behavior.}
Throughout this section, let $F:X\dashrightarrow Y$ be a dominant
meromorphic map between compact K\"ahler surfaces. Following~\cite[\S
34.7]{manin}, we introduce the action of $F$ on Weil and Cartier
classes. We then describe the continuity properties of these actions
on the Hilbert space $\Ltwo(\fX)$.

For each blow-up $Y_\varpi$ of $Y$, there
exists a blow-up $X_\pi$ of $X$ such that the induced map 
$X_\pi\to Y_\varpi$ is holomorphic. The associated push-forward
$\hodge(X_\pi)\to\hodge(Y_\varpi)$ and pull-back
$\hodge(Y_\varpi)\to\hodge(X_\pi)$ are compatible with the
projective and injective systems defined by push-forwards and
pull-backs that define Weil and Cartier classes respectively, so we
can consider the induced morphisms on the respective 
projective and inductive limits.
\begin{defi} 
  Given $F:X\dashrightarrow Y$ as above, 
  we denote by $F_*:W(\fX)\to W(\fY)$ the induced push-forward operator,
  and by $F^*:C(\fY)\to C(\fX)$ the induced pull-back operator.
\end{defi}
Concretely, if $\a\in W(\fX)$ is a Weil class, the
incarnation of $F_*\a\in W(\fY)$ on a given blow-up $Y_\varpi$ is the
push-forward of $\a_\pi\in\hodge(X_\pi)$ by the induced map
$X_\pi\to Y_\varpi$ for any $\pi$ such that the latter map
is holomorphic. Similarly, if $\b\in C(\fY)$ is
a Cartier class determined on a blow-up $Y_\varpi$, its pull-back
$F^*\b\in C(\fX)$ is the Cartier class determined on $X_\pi$ by the
pull-back of $\b_\varpi\in\hodge(Y_\varpi)$ 
by the induced map $X_\pi\to Y_\varpi$, 
whenever the latter is holomorphic.

These constructions are functorial, \ie $(F\circ G)_*=F_*\circ G_*$
and $(F\circ G)^*=G^*\circ F^*$, and compatible with the duality
between $C$ and $W$, since this is true for each holomorphic map
$X_\pi\to Y_\varpi$. In other words, for any $\a\in W(\fX)$ and
$\b\in C(\fY)$, we have $(F_*\a\cdot\b)=(\a\cdot F^*\b)$.

We also see that $F_*$ preserves nef
and psef Weil classes, and that $F^*$ preserves nef and psef Cartier
classes. Indeed, the pull-back and push-forward by a surjective
holomorphic map both preserve nef and psef $(1,1)$-classes in
dimension two.
\begin{rmk}\label{rem:usualpp}
  If $\pi:X_\pi\to X$ and $\varpi:Y_\varpi\to Y$ are arbitrary
  blowups, then the pullback operator
  $\hodge(Y_\varpi)\to\hodge(X_\pi)$ usually 
  associated to the meromorphic map $X_\pi\dashrightarrow Y_\varpi$
  is given by the restriction of $F^*:C(\fY)\to C(\fX)$ to
  $\hodge(Y_\varpi)$, followed by the projection 
  of $C(\fX)$ onto $\hodge(X_\pi)$. 
  Similarly, the pushforward operator
  $\hodge(X_\pi)\to\hodge(Y_\varpi)$ usually 
  associated to $X_\pi\dashrightarrow Y_\varpi$
  is given by the restriction of $F_*:W(\fX)\to W(\fY)$ to
  $\hodge(X_\pi)$, followed by the projection 
  of $W(\fY)$ onto $\hodge(Y_\varpi)$. 
\end{rmk}

The intersection forms on $C(\fX)$ and $C(\fY)$ are related by $F^*$
as follows: $(F^*\b^2)=e(F)(\b^2)$, where $e(F)>0$ is the topological
degree of $F$. In view of the universal property of completions
mentioned in \S\ref{L2} on p.\pageref{blob}, we get
\begin{prop} 
  The pull-back $F^*:C(\fY)\to C(\fX)$ extends to a continuous
  operator $F^*:\Ltwo(\fY)\to\Ltwo(\fX)$ such that
  $((F^*\b)^2)=e(F)(\b^2)$ for each $\b\in\Ltwo(\fY)$. 
  By duality, the push-forward $F_*:W(\fX)\to W(\fY)$ 
  induces a continuous operator
  $F_*:\Ltwo(\fX)\to\Ltwo(\fY)$, such that 
  $(F_*\a\cdot\b)=(\a\cdot F^*\b)$ for any
  $\a,\b\in\Ltwo(\fX)$.
\end{prop}
Next we show that the pull-back $F^*:C(\fY)\to C(\fX)$
continuously extends to Weil classes 
and---dually---the push-forward 
$F_*:W(\fX)\to W(\fY)$ preserves Cartier classes.

In doing so, we shall repeatedly use 
a consequence of the result of Zariski 
already mentioned before. Namely, given $F:X\dashrightarrow Y$ and
a blowup $\pi:X_\pi\to X$, there
exists a blow-up $Y_\varpi$ of $Y$ such that the induced meromorphic
map $X_\pi\dashrightarrow Y_\varpi$ 
does not contract any curve to a point. 
\begin{lem} 
  Suppose $\pi:X_\pi\to X$, and
  $\varpi: Y_\varpi \to Y$ are two blow-ups such that the induced
  meromorphic map $X_\pi\dashrightarrow Y_\varpi$ does not contract
  any curve to a point.  Then for each Cartier class $\b\in C(\fY)$, 
  the incarnations of $F^*\b$ and $F^*\b_{\varpi}$ on $X_\pi$
  coincide.
\end{lem}
\begin{proof}
  Any Cartier class is a difference of nef Cartier classes
  so we may assume $\beta$ is nef and determined in some blowup
  $\varpi'$ dominating $\varpi$. Pick $\pi'$ dominating 
  $\pi$ such that the induced map $X_{\pi'}\to Y_{\varpi'}$
  is holomorphic. 
  Set $\a\=F^*(\beta_\varpi-\beta)$. 
  Then $\a\in C(\fX)$ is psef and determined in $X_{\pi'}$.
  We must show that $\a_\pi=0$.
  If $\a_\pi\ne0$, then $\a\ge\lambda C$, 
  where $\lambda>0$ and $C$ is the class of an irreducible
  curve on $X_\pi$. Now $C$ is not contracted by
  $X_\pi\dashrightarrow Y_\varpi$ so the incarnation of
  $F_*\a$ on $Y_\varpi$ is nonzero.
  But this is a contradiction, since this incarnation
  equals $e(F)(\b_\varpi-\b)_\varpi=0$.
\end{proof}
\begin{cor}\label{c:pull} 
  The pull-back operator
  $F^*: C(\fY) \to C(\fX)$ continuously extends to 
  $F^*:W(\fY)\to W(\fX)$, and preserves
  nef and psef Weil classes.
  
  More precisely, if $X_\pi$ is a given blow-up of $X$, and $Y_{\varpi}$
  is a blow-up of $Y$ such that the induced meromorphic map 
  $X_\pi\dashrightarrow Y_{\varpi}$ does not contract curves, 
  then for any Weil class $\gamma\in W(\fY)$, 
  one has $(F^*\gamma)_\pi=(F^*\gamma_\varpi)_\pi$.
\end{cor}
\begin{cor}\label{c:pushcart} 
  The push-forward operator
  $F_*:W(\fX)\to W(\fY)$ preserves Cartier classes. More precisely, if
  $\a\in C(\fX)$ is a Cartier class determined on some $X_\pi$, then
  $F_*\a$ is Cartier, determined on $Y_\varpi$ as soon as the induced
  meromorphic map $X_\pi\dashrightarrow Y_{\varpi}$ does not contract
  curves.
\end{cor} 
\begin{proof}
  For any $\b\in C(\fY)$, the incarnations of
  $F^*\b$ and $F^*\b_{\varpi}$ on $X_{\pi}$ coincide
  by Corollary~\ref{c:pull}.  
  Hence 
  \begin{equation*}
    (F_*\a\cdot\b)
    =(\a\cdot F^*\b)
    =(\a\cdot F^*\b_{\varpi})
    =(F_*\a\cdot\b_{\varpi})
    =((F_*\a)_\varpi\cdot\b).
  \end{equation*} 
  As this holds for any Cartier class $\b\in C(\fY)$ we must have
  $F_*\a=(F_*\a)_\varpi$ by Proposition~\ref{p:dual}.
\end{proof}
%
%
%%%%%%%%%%%%%%%%%%%%%%%%%%%%%%%%%%%%%%%%%%%%%%%%%%%%%%%%%%%%%%%%%%%
%
%
\section{Dynamics}
Now consider a dominant meromorphic self-map $F:X\dashrightarrow X$ of
a compact K{\"a}hler surface $X$.  Write $\la_2=e(F)$ for the
topological degree of $F$.  If $\om\in\Nef(\fX)$ is a nef Weil class
such that $(\omega^2)>0$, we define the degree of $F$ with respect to
$\omega$ as
\begin{equation*}
  \deg_\om(F)\=(F^*\om\cdot \om)=(\om\cdot F_*\om).
\end{equation*}
This coincides
with the usual notion of degree when $X=\P^2$ and $\omega$ is the
Cartier class determined by a line on $\P^2$.

\begin{prop}\label{prop:dinfty}
  The limit
  \begin{equation}\label{eq:dyndeg}
    \la_1\=\la_1(F)\=\lim_{n\to\infty}\deg_{\om}(F^n)^{\frac1n}
  \end{equation}
  exists and does not depend on the choice of the nef class
  $\om\in\Nef(\fX)$ with $(\om^2)>0$. Moreover, $\la_1$ is invariant
  under bimeromorphic conjugacy and $\la_1^2\ge\la_2$.
\end{prop}
The result above is well known but we include the proof for
completeness.  We call $\la_1$ the \emph{asymptotic degree} of $F$.
It is also known as the \emph{first dynamical degree} and can be
computed (see~\cite{DF}) as $\la_1=\lim_{n\to\infty}\rho_n^{1/n}$,
where $\rho_n$ is the spectral radius of $F^n$
acting on $\hodge(X)$ by pullback or push-forward 
(\cf~Remark~\ref{rem:usualpp}).
\begin{proof}[Proof of Proposition~\ref{prop:dinfty}]
  Upon scaling $\om$, we can assume that
  $(\om^2)=1$. By~\eqref{eq:nefbounds} we then have
  $G^*\om\le2(G^*\om\cdot \om) \, \om$ for any dominant mapping
  $G:X\dashrightarrow X$.  Applying this with $G=F^m$ yields
  \begin{equation*}
    \deg_{\om}{F^{n+m}}
    =(F^{n*}F^{m*}{\om}\cdot {\om})
    \le2(F^{n*}{\om}\cdot {\om})(F^{m*}{\om}\cdot {\om})
    =2\deg_{\om}(F^n)\deg_{\om}(F^m)
  \end{equation*}
  This implies (see \eg~\cite[Prop.~9.6.4]{KH}) that the limit
  in~\eqref{eq:dyndeg} exists. Let us temporarily denote it by
  $\lambda_1(\omega)$.  If $\om'\in C(\fX)$ is another nef class with
  $(\om'^2)>0$, then it follows from~\eqref{eq:nefbounds} that
  $\om'\le C\om$ for some $C>0$.  By Proposition~\ref{p:monotone},
  this gives
  \begin{equation*}
    \deg_{\om'}F^n
    =(F^{n*}\om'\cdot\om')
    \le C^2(F^{n*}\om\cdot\om)
    =C^2\deg_{\om}F^n
     \end{equation*}
  Taking $n$th roots and letting $n\to\infty$ shows that
  $\lambda_1(\om')\le\lambda_1(\om)$, and thus
  $\lambda_1(\om')=\lambda_1(\om)$ by symmetry, so that $\la_1$ is
  indeed independent of ${\om}$.  It is then invariant by
  bimeromorphic conjugacy, since $\fX$ and all the spaces attached to
  it are.
  
  Finally, Proposition~\ref{p:monotone} 
  yields $F^{n*}\om\leq 2(F^{*n}\om\cdot\om)\,\om$, which implies
  \begin{equation*}
    e(F)^n
    =e(F^n)
    =(F^{n*}\om^2)
    \leq 4(F^{n*}\om\cdot\om)^2
    =4\deg_\om(F^n)^2
  \end{equation*}
  and letting $n\to\infty$ yields $\lambda_2=e(F)\leq\lambda_1^2$.
\end{proof}
%
%%%%%%%%%%%%%%%%%%%%%%%%%%%%%%%%%%%%%%%%%%%%%%%%%%%%%%%%%%%%%%%%%%%
%
\subsection{Existence of eigenclasses}
To begin with we do not assume $\la_1^2>\la_2$.
\begin{thm}\label{thm:eigenclass}
  Let $F:X\dashrightarrow X$ be any dominant meromorphic selfmap of a
  smooth K{\"a}hler surface $X$ with asymptotic degree $\la_1$.  Then
  we can find nonzero nef Weil classes $\theta_\ast$ and $\theta^\ast$
  with $F_*\theta_\ast=\la_1\theta_\ast$ and
  $F^*\theta^\ast=\la_1\theta^\ast$.
\end{thm}
Note that by Proposition~\ref{p:monotone}, both classes $\theta_\ast,
\theta^\ast$ belong to $\Ltwo(\fX)$.
\begin{proof}
  We shall use the push-forward and pull-back operators
  \begin{equation*}
    S_\pi:\hodge(X_\pi)\to\hodge(X_\pi)
    \quad\text{and}\quad
    T_\pi:\hodge(X_\pi)\to\hodge(X_\pi)
  \end{equation*}      
  usually associated to the meromorphic map 
  $X_\pi\dashrightarrow X_\pi$ induced by $F$
  for a given blowup $\pi:X_\pi\to X$.
  Thus $S_\pi$ (resp.\ $T_\pi$) is
  the restriction to $\hodge(X_\pi)$
  of $F_*:C(\fX)\to C(\fX)$ 
  (resp.\ $F^*:C(\fX)\to C(\fX)$) 
  followed by the projection 
  $C(\fX)\to\hodge(X_\pi)$, cf.\ Remark~\ref{rem:usualpp}.
  These operators are typically denoted $F_*$ and $F^*$ 
  in the literature, but here that notation would conflict with 
  the corresponding operators on $C(\fX)$ or $W(\fX)$.

  The spectral radius
  $\rho_\pi>0$ of $T_\pi$ can be computed as follows: if
  $\theta\in\hodge(X_\pi)$ is any nef class with $(\theta^2)>0$, then
  $(T_\pi^n\theta\cdot\theta)^{1/n}\to\rho_\pi$ as $n\to\infty$.
  \begin{lem} 
    We have $\lambda_1\le\rho_{\pi'}\le\rho_{\pi}$ 
    for all $\pi'\ge\pi$.
  \end{lem}
  \begin{proof} 
    Let $\theta\in C(\fX)$ be a given nef 
    class determined on $X_{\pi'}$ with $(\theta^2)>0$, so that
    $\theta\le\theta_{\pi}$ by Proposition~\ref{p:incr}. Then
    $T_{\pi'}\theta$ is the incarnation on $X_{\pi'}$ of the nef class
    $F^*\theta$ on $X_{\pi'}$, and 
    $T_{\pi}\theta_{\pi}$ is the incarnation
    on $X_{\pi}$ of the nef class $F^*\theta_{\pi}\ge F^*\theta$, 
    thus $F^*\theta\le T_{\pi'}\theta\le T_{\pi}\theta_{\pi}$ holds by
    Proposition~\ref{p:incr}. By induction we get 
    $F^{n*}\theta\le T_{\pi'}^n\theta\le T_{\pi}^n\theta_{\pi}$ 
    for all $n$, hence
    $(F^{n*}\theta\cdot\theta)^{1/n}\le(T_{\pi'}^n\theta\cdot\theta)^{1/n}
    \le(T_{\pi}^n\theta_{\pi}\cdot\theta_{\pi})^{1/n}$ by
    Proposition~\ref{p:monotone}, 
    and $\lambda_1\le\rho_{\pi'}\le\rho_{\pi}$
    follows by letting $n\to\infty$.
  \end{proof}
  
  Now the set of nef classes in $\hodge(X_\pi)$ is a closed convex cone
  with compact basis invariant by $T_\pi$, thus a Perron-Frobenius
  type argument (see~\cite[Lemma~1.12]{DF}) establishes the existence
  of a non-zero nef class $\theta(\pi)\in\hodge(X_\pi)$ with
  $T_\pi\theta(\pi)=\rho_\pi\theta(\pi)$.
  
  If we identify $\theta(\pi)$ with the nef Cartier class it
  determines, this says that the nef Cartier classes $F^*\theta(\pi)$
  and $\rho_\pi\theta(\pi)$ have the same incarnation on $X_\pi$. We
  have thus obtained approximate eigenclasses, and the plan is now
  to get the desired class $\theta^*$ as a limit of classes of the
  form $\theta(\pi)$. We will then explain how to modify the argument
  to construct $\theta_\ast$.
  
  We normalize $\theta(\pi)$ by $(\theta(\pi)\cdot \om)=1$ for a fixed
  class $\om\in C(\fX)$ determined by a K\"ahler class on $X$ with
  $(\om^2)=1$, so that the $\theta(\pi)$ all lie in a compact subset of
  the nef cone $\Nef(\fX)$ by Proposition~\ref{p:compact}.
  
  Let $X=X_0\leftarrow X_1\leftarrow\dots$ be an infinite sequence of
  blow-ups, such that the lift of $F$ 
  as a map from $X_{n+1}$ to $X_n$ is holomorphic for $n\ge0$. 
  
  For each $n$, let $\rho_n$ denote the spectral radius of $T_n$ on
  $\hodge(X_n)$ as above, and pick a non-zero nef Cartier class
  $\theta_n\in C(\fX)$ determined on $X_n$ and such that $T_n \theta_n =
  \rho_n \theta_n$. Then $F^* \theta_n$ is a
  Cartier class determined in $X_{n+1}$, and by definition $T_n
  \theta_n$ is the incarnation of this class in $X_n$. Therefore 
  $F^*\theta_n$ and $\rho_n\theta_n$ coincide on $X_n$. By
  Proposition~\ref{p:conv}, it follows that $F^*\theta_n-\rho_n\theta_n$
  converges to $0$ in $W(\fX)$ as $n\to\infty$.
  
  We have seen above that $\rho_n$ is a decreasing sequence. Let
  $\rho_\infty\=\lim\rho_n$, so that $\rho_\infty\ge\lambda_1$ by the
  above lemma. Since the $\theta_n$ lie in a compact subset of
  $\Nef(\fX)$, we can find a cluster point $\theta^\ast$ for the
  sequence $\theta_n$, which is also a nef Weil class with
  $(\theta^{\ast}\cdot\omega)=1$. Since $F^*\theta_n-\rho_n\theta_n$
  converges to $0$ in $W(\fX)$, it follows that
  $F^*\theta^{\ast}=\rho_\infty\theta^{\ast}$.
  
  To complete the proof we will show that $\rho_\infty=\la_1$. In fact,
  if $\a\in W(\fX)$ is any non-zero nef eigenclass of $F^*$ with
  $F^*\a=t\,\a$ for some $t\ge 0$, then $t\le\la_1$. Indeed, we have
  $\a\le C\om$ for some $C>0$ by Proposition~\ref{p:nefbounds}, and
  it follows that $(F^{n*}\om\cdot {\om}) \ge C^{-1}(F^{n*}\a\cdot
  {\om}) =C^{-1}t^n(\a\cdot\om)$.  Taking $n$th roots and letting
  $n\to\infty$ yields $\la_1\ge t$, as was to be shown.
  
  \smallskip
  In order to construct $\theta_*$, we modify the above argument as
  follows.
  Let $S_\pi:\hodge(X_\pi)\to\hodge(X_\pi)$ 
  be the push-forward operator defined above.
  As $F^*$ and $F_*$ are adjoint to each other with respect to the
  intersection pairing, it follows that $S_\pi$ and $T_\pi$ are adjoint
  with respect to Poincar\'e duality on $\hodge(X_\pi)$, so that they
  have the same spectral radius $\rho_\pi$. By Perron-Frobenius, there
  exists a non-zero nef class $\vartheta(\pi)\in\hodge(X_\pi)$ such that
  $S_\pi\vartheta(\pi)=\rho_\pi\vartheta(\pi)$.
  
  Now pick $X=X_0\leftarrow X_1\leftarrow\dots$ an infinite sequence
  of blow-ups such that the lifts of $F$ from $X_n$ to $X_{n+1}$ do
  not contract any curves.  For each $n$, we get a nef class
  $\vartheta_n\in C(\fX)$ determined on $X_n$ normalized by
  $(\vartheta_n\cdot\om)=1$.  By Corollary~\ref{c:pushcart}, the class
  $F_*\vartheta_n$ is determined in $X_{n+1}$, so $F_*\vartheta_n$ and
  $\rho_n\vartheta_n$ coincide in $X_n$.  Proposition~\ref{p:conv}
  then shows that $F_*\vartheta_n-\rho_n\vartheta_n$ converges to $0$
  in $W(\fX)$ as $n\to\infty$, hence $\theta_\ast\in\Nef(\fX)$ can be
  taken to be any cluster value of $\vartheta_n$.
\end{proof}
\begin{rmk}
  When $K_X$ is not psef (\ie if $X$ is rational or ruled) we may also
  achieve $(\theta_\ast\cdot K_{\fX})\le0$. To see this, first note
  that $F^*K_{\fX}\le K_{\fX}$ as classes in $W(\fX)$, since
  $K_{X_{\pi'}}-F^*K_{X_\pi}$ is represented by the effective 
  zero-divisor of
  the Jacobian determinant of the map $X_{\pi'}\to X_\pi$ induced
  by $F$ assuming this is holomorphic.
  Now for each blow-up $X_\pi$, let $C_\pi$ be the set of
  nef classes $\a\in\hodge(X_\pi)$ such that $(\a\cdot K_{\fX})\le0$.
  Then $C_\pi$ is a closed convex cone with compact basis, and is not
  reduced to $0$ since $K_\fX$ is not psef. It is furthermore
  invariant by $S_\pi$. Indeed, if $\a\in\hodge(X_\pi)$ is a nef
  class, we have
  \begin{equation*}
    (S_\pi\a\cdot K_{\fX}) 
    =(F_*\a\cdot K_{X_\pi}) 
    \le(F_*\a\cdot K_{\fX}) 
    =(\a\cdot F^*K_{\fX})
    \le(\a\cdot K_{\fX}).
  \end{equation*}
  We can thus assume
  that the non-zero eigenclasses $\vartheta_n$ in the proof above belong to
  $C_n$, and we get $(\theta_\ast\cdot K_{\fX})\leq 0$.
  
  The same argument does \emph{not} work for $\theta^\ast$, since
  $F_*K_\fX\le K_{\fX}$ does not hold in general.
\end{rmk}
%
%%%%%%%%%%%%%%%%%%%%%%%%%%%%%%%%%%%%%%%%%%%%%%%%%%%%%%%%%%%%%%%%%%%
%
\subsection{Spectral properties}
Theorem~\ref{thm:eigenclass} asserts the existence of eigen\-clas\-ses
for $F_*$ and $F^*$ with eigenvalue $\la_1$.  We now further analyze
the spectral properties under the assumption that $\la_1^2>\la_2$.

\begin{thm}\label{thm:spectral}
  Assume $\la_1^2>\la_2$. Then the non-zero nef Weil classes
  $\theta_\ast, \theta^\ast\in \Ltwo(\fX)$ such that $F^*\theta^\ast=\la_1
  \theta^\ast$ and $F_*\theta_\ast=\la_1\theta_\ast$ are unique up to
  scaling. We have $(\theta_\ast\cdot\theta^\ast)>0$ and
  $(\theta^{\ast 2})=0$. We rescale them so that
  $(\theta_\ast\cdot\theta^\ast)=1$. Let $\cH\subset\Ltwo(\fX)$ be the
  orthogonal complement of $\theta^*$ and $\theta_*$, so that we have
  the decomposition
  $\Ltwo(\fX)=\R\theta^\ast\oplus\R\theta_\ast\oplus\cH$. The
  intersection form is negative definite on $\cH$, and
  $\|\a\|^2:=-(\a^2)$ defines a Hilbert norm on $\cH$.  The actions of
  $F^*$ and $F_*$ with respect to this decomposition are as follows:
  \begin{itemize}
  \item[(i)]
    The subspace $\cH$ is $F^*$-invariant and 
    \begin{equation*}
      \begin{cases}
        F^{n*}\theta^\ast
        =\la_1^n\,\theta ^\ast;\\
        F^{n*}\theta_\ast
        =(\frac{\la_2}{\la_1})^n\theta_\ast
        +(\theta_*^2)\,\la_1^n(1-(\frac{\la_2}{\la_1^2})^n)\,\theta^\ast
        +h_n\\ \hfill \text{with}\  h_n\in\cH,\ 
        \|h_n\|=O(\la_2^{n/2});\\
        \|F^{n*}h\|=\la_2^{n/2}\|h\|\ 
        \text{for all $h\in\cH$}.
      \end{cases}
    \end{equation*}
  \item[(ii)]
    The subspace $\cH$ is not $F_*$-invariant in general, but
    \begin{equation*}
      \begin{cases}
        F^n_*\theta _\ast=\la_1^n\theta_\ast;\\
        F^n_*\theta^\ast= (\frac{\la_2}{\la_1})^n\theta^\ast;\\
        \|F^n_*h\|\le C\la_2^{n/2}\|h\|\ 
        \text{for some $C>0$ and all $h\in\cH$}.
      \end{cases}
    \end{equation*}
  \end{itemize}
\end{thm}

\begin{cor}\label{cor:asympt}
  For any Weil class $\a\in\Ltwo(\fX)$, we have 
  \begin{equation*}
    \frac1{\la_1^n}F^{n*}\a
    =(\a\cdot\theta_\ast)\theta^\ast
    +O((\frac{\la_2}{\la_1^2})^{n/2})
    \text{ and }
    \frac1{\la_1^n}F^n_*\a
    =(\a\cdot\theta^*)\theta_\ast
    +O((\frac{\la_2}{\la_1^2})^{n/2}).
  \end{equation*}
\end{cor}
\begin{proof}
  The decomposition of $\a$ in 
  $\Ltwo(\fX)=\R\theta^\ast\oplus\R\theta_\ast\oplus\cH$
  is given by
  \begin{equation}\label{eq:decomp}
    \a
    =((\a\cdot\theta_\ast)-(\a\cdot\theta^\ast)(\theta_*^2))\theta^\ast
    +(\a\cdot\theta^\ast)\theta_\ast
    +\a_0,
  \end{equation}
  where $\a_0\in\cH$. The result follows from~\eqref{eq:decomp}
  using~(i) and~(ii) above.
\end{proof}
\begin{proof}[Proof of the Main Theorem]
  Applying Corollary~\ref{cor:asympt} to
  $\a={\om}$ (which is nef, hence in $\Ltwo(\fX)$) 
  gives
  \begin{equation*}
    \deg_{\om}(F^n)
    =(F^{n*}{\om}\cdot {\om})
    =(\omega\cdot\theta^*)(\omega\cdot\theta_*)\la_1^n+O(\la_2^{n/2}),
  \end{equation*}
  This completes the proof with $b:=(\om\cdot\theta^*)(\om\cdot\theta_*)$. 
\end{proof}
\begin{proof}[Proof of Theorem~\ref{thm:spectral}]
  Using Theorem~\ref{thm:eigenclass}, 
  we may find nonzero nef Weil classes
  $\theta_\ast, \theta^\ast$ such that 
  $F_*\theta_\ast=\la_1\theta_\ast$ 
  and $F^*\theta^\ast=\la_1\theta^\ast$.
  Fix two such classes for the duration of the proof.
  In the end we shall see that they are unique up to scaling.

The proof amounts to a series of simple arguments using general facts
for transformations of a complete vector space endowed with a Minkowski
form. We provide the details for the benefit of the reader.

First note that
  $\lambda_1F_*\theta^*=F_*F^*\theta^*=\lambda_2\theta^*$, so that
  $F_*\theta^*=(\lambda_2/\lambda_1)\theta^*$. Since
  $F_*\theta_*=\lambda_1\theta_*$ and $\lambda_1^2>\lambda_2$, it
  follows that $\theta^*$ and $\theta_*$ cannot be proportional.

  Applying the relation $(F^*\a^2)=\lambda_2(\a^2)$ to $\a=\theta^*$
  yields $\lambda_1^2(\theta^{*2})=\lambda_2(\theta^{*2})$, and thus
  $(\theta^{*2})=0$ since $\la_1^2>\la_2$. By the Hodge index theorem,
  $\theta_*$ and $\theta^*$ would thus have to be proportional if they
  were orthogonal. We infer that $(\theta^*\cdot\theta_*)>0$, and we
  rescale $\theta^*$ so that $(\theta^*\cdot\theta_*)=1$.
  
   Let us first prove the properties in~(i) for the pullback.
  As both $\theta_\ast$ and $\theta^\ast$ are eigenvectors for
  $F_*$, the space $\cH$ is invariant under $F^*$.
  Using~\eqref{eq:decomp} and the invariance properties of $\theta_\ast$ and 
  $\theta^\ast$, we get
  \begin{equation}\label{eq:pull+}
    F^*\theta_\ast
    =\frac{\la_2}{\la_1}\theta_\ast
    +\la_1(1-\frac{\la_2}{\la_1^2})(\theta_*^2)\theta^\ast
    +h_1,
  \end{equation}
  where $h_1\in\cH$. Inductively,~\eqref{eq:pull+} gives
  \begin{equation}\label{eq:pull2}
    F^{n*}\theta_\ast
    =(\frac{\la_2}{\la_1})^n\theta_\ast
    +\la_1^n(1-(\frac{\la_2}{\la_1^2})^n)(\theta_*^2)\theta^\ast
    +h_n,
  \end{equation}
  where $h_{n+1}=F^*h_n+(\la_2/\la_1)^nh_1 \in \cH$. Using that
  $\|F^*h\|^2=\lambda_2\|h\|^2$ on $\cH$, we get
  $\|h_{n+1}\|\leq\la_2^{1/2}\|h_n\|+(\la_2/\la_1)^n\|h_1\|$, which is
  easily seen to imply $\|h_n\|=O(\la_2^{n/2})$ since
  $\sum_k(\lambda_2^{1/2}/\lambda_1)^k<+\infty$.  This concludes the
  proof of~(i).
  
  Let us now turn to the push-forward operator. The first two
  equations are clear.  As $\theta_\ast$ may not be an eigenvector for
  $F^*$, $\cH$ need not be invariant by $F_*$, but since $F_*h$ is
  orthogonal to $\theta^*$ for any $h\in\cH$, we can write
  $F^n_*h=a_n\theta^*+g_n$, with $a_n=(F^{n*}\theta_*\cdot h)$ and
  $g_n\in\cH$. We have seen that $F^{n*}\theta_*=h_n$ modulo
  $\theta^*,\theta_*$ with $\|h_n\|=O(\la_2^{n/2})$, thus
  $|a_n|=|(h_n\cdot h)|\le C\la_2^{n/2}\|h\|$. On the other hand, we
  have $(g_n^2)=(F^{n*}g_n\cdot h)$, and thus $\|g_n\|^2\leq
  \la_2^{n/2}\|g_n\|\|h\|$, and this shows that 
  $\|F^n_*h\|\le C\la_2^{n/2}\|h\|$.
\end{proof}

\begin{rmk}\label{rem:bomega}
  It follows from the proof of the Main Theorem that there exist
  nef classes $\a_\ast, \a^\ast\in\hodge(X)$ such that
  for any K\"ahler classes $\om$, $\om'$ on $X$, we have 
  \begin{equation*}
    \frac{\deg_{\om}(F^n)}{\deg_{\om'}(F^n)}
    =\frac{(\alpha^\ast\cdot\om)_X\, (\alpha_\ast\cdot\om)_X}
    {(\alpha^\ast\cdot\om')_X\, (\alpha_\ast\cdot\om')_X}
    +O((\frac{\la_2}{\la_1^2})^{n/2}).
%    +O\left(\left(\frac{\la_2}{\la_1^2}\right)^{n/2}\right).
  \end{equation*}
  Indeed, we can take
  $\a^\ast$ and $\a_\ast$ as the incarnations in $X$ of  
  $\theta^\ast$ and $\theta_\ast$, respectively.
\end{rmk}
\begin{rmk}
  When $F$ is bimeromorphic we have  
  $\theta_\ast(F)=\theta^\ast(F^{-1})$,
  hence $(\theta_\ast^2)=0$.
  However in general we may have 
  $(\theta_\ast^2)>0$.
  For example, let $F$ be any polynomial map of $\C^2$ 
  whose extension to $\P^2$ is not holomorphic but 
  does not contract any curve. 
  If $\omega$ is the class of a line on $\P^2$, then
  $\deg_\omega(F)>\sqrt{\la_2}>1$. On the other hand,
  $F_*\omega=\deg_\omega(F)\omega$ 
  by Corollary~\ref{c:pushcart},
  so $\la_1=\deg_\omega(F)$, $\theta_\ast=\omega$ and 
  $(\theta_\ast^2)=1$.
\end{rmk}
\begin{rmk}
  The case when $\theta_\ast$ (or $\theta^\ast$) is Cartier is very
  special. For example, when $F$ is bimeromorphic, it follows
  from~\cite[Theorem~0.4]{DF} that $\theta_\ast$ (or, equivalently,
  $\theta^*$) is Cartier iff $F$ is \emph{biholomorphic} in some
  birational model. In the general non-invertible case, similar
  rigidity results are expected, see~\cite{cantat2} for work in this
  direction.

  Note also that $F$ being algebraically stable in some birational
  model does not imply that the eigenclasses are Cartier. We do not
  know whether having a Cartier eigenclass implies algebraic stability
  in some model, but having a Cartier eigenclass has many of the same
  consequences as stability: $\la_1$ is an algebraic integer and the
  sequence of degrees $(\deg_\om F^n)_1^\infty$ satisfies a linear
  recurrence relation.
\end{rmk}
%
%
%%%%%%%%%%%%%%%%%%%%%%%%%%%%%%%%%%%%%%%%%%%%%%%%%%%%%%%%%%%%%%%%%%%
%
%

\end{document}